\documentclass{article}
\usepackage{amssymb,amsmath,stmaryrd,array,graphicx,mathtools}
\title{A Theory of Particular Sets}
\author{Paul Blain Levy, University of Birmingham}
\usepackage[margin=2cm]{geometry}
\usepackage[matrix,arrow,curve,rotate]{xy}
\newcommand{\condsing}[2]{\setbr{{#1}\mid\mid{#2}}}

\newcommand{\interline}{\\ \vspace{-2ex} \\}

\newcommand{\pset}{\mathcal{P}}
\newcommand{\eqdef}{\stackrel{\mbox{\rm {\tiny def}}}{=}}
\newcommand{\cata}{\ensuremath{\mathcal{A}}}
\newenvironment{spaceout}[1]{\begin{displaymath}\setlength{\extrarowheight}{3pt}\begin{array}{#1}}{\end{array}\setlength{\extrarowheight}{0pt}\end{displaymath} \noindent}

\newcommand{\nonsenseconv}[2]{{#1}\; \rightarrowtriangle_{\mathsf{NC}}\; {#2}}

\newcommand{\betwixt }{\hspace{1em}}

\newcommand{\textb}[1]{\text{\textbf{#1}}}
\newcommand{\nats}{\mathbb{N}}
\newcommand{\bd}{\ }
\newcommand{\setbr}[1]{\{{#1}\}}
\newcommand{\tuple}[1]{\langle {#1} \rangle}

\newcommand{\zf}{ZFC}

\newcommand{\rotiota}{\rotatebox[origin=c]{180}{$\iota$}}
\newcommand{\totall}{\mathfrak{T}}
\newcommand{\allsets}{\mathfrak{S}}
\newcommand{\vars}{\mathsf{Vars}}
\newcommand{\isset}[1]{\mathsf{IsSet}({#1})}

\newcommand{\Intro}{\mathcal{I}}
\newcommand{\Elim}{\mathcal{E}}
\newcommand{\parenth}[1]{({#1})}
\newcommand{\textp}[1]{\text{({#1})}}

\newcommand{\uniqel}[1]{\mathsf{UE}({#1})}

\newcommand{\itset}[2]{\mathsf{IR}({#1}/{#2})} 
\newcommand{\wfint}[3]{\mathsf{W}_{\in\! {#2}}^{{#3}}({#1})}

\newcommand{\prl}[4]{\mathsf{WR}_{\in \!{#1}}^{{#2}/ {#3}}({#4})}

\newcommand{\ifflong}{\ \Leftrightarrow \ }

\newtheorem{definition}{Definition}

\newcommand{\dvar}[1]{\mathsf{Decl}({#1})}
\newcommand{\dvarg}[1]{\dvar{#1}}

\newcommand{\lco}{,\,}
\newcommand{\llco}{,\,}
\newcommand{\hyps}[1]{\mathsf{Hyp}({#1})}

\newcommand{\colin}{\!:\!}

\newcommand{\smin}{\!\in\!}

\newcommand{\tops}{TOPS}

\newcommand{\bbot}{\emptyset}
\newcommand{\False}{\mathsf{False}}
\newcommand{\True}{\mathsf{True}}

\newcommand{\vdasht}{\vdash^{\mathsf{t}}}
\newcommand{\vdashp}{\vdash^{\mathsf{p}}}
\newcommand{\sht}[2]{{#1} \vdasht {#2}}
\newcommand{\shp}[2]{{#1} \vdashp {#2}}

\newcommand{\bnfgo}{::=\hspace{0.7em}}
\newcommand{\bnf}{\hspace{0.5em} | \hspace{0.5em}}

\newcommand{\existsu}{\exists!}

\newdimen\proofrulebreadth \proofrulebreadth=.05em
\newdimen\proofdotseparation \proofdotseparation=1.25ex
\newdimen\proofrulebaseline \proofrulebaseline=2ex
\newcount\proofdotnumber \proofdotnumber=3
\let\then\relax
\def\hfi{\hskip0pt plus.0001fil}
\mathchardef\squigto="3A3B
\newif\ifinsideprooftree\insideprooftreefalse
\newif\ifonleftofproofrule\onleftofproofrulefalse
\newif\ifproofdots\proofdotsfalse
\newif\ifdoubleproof\doubleprooffalse
\let\wereinproofbit\relax
\newdimen\shortenproofleft
\newdimen\shortenproofright
\newdimen\proofbelowshift
\newbox\proofabove
\newbox\proofbelow
\newbox\proofrulename
\def\shiftproofbelow{\let\next\relax\afterassignment\setshiftproofbelow\dimen0 }
\def\shiftproofbelowneg{\def\next{\multiply\dimen0 by-1 }%
\afterassignment\setshiftproofbelow\dimen0 }
\def\setshiftproofbelow{\next\proofbelowshift=\dimen0 }
\def\setproofrulebreadth{\proofrulebreadth}

\def\prooftree{
\ifnum  \lastpenalty=1
\then   \unpenalty
\else   \onleftofproofrulefalse
\fi
\ifonleftofproofrule
\else   \ifinsideprooftree
        \then   \hskip.5em plus1fil
        \fi
\fi
\bgroup
\setbox\proofbelow=\hbox{}\setbox\proofrulename=\hbox{}%
\let\justifies\proofover\let\leadsto\proofoverdots\let\Justifies\proofoverdbl
\let\using\proofusing\let\[\prooftree
\ifinsideprooftree\let\]\endprooftree\fi
\proofdotsfalse\doubleprooffalse
\let\thickness\setproofrulebreadth
\let\shiftright\shiftproofbelow \let\shift\shiftproofbelow
\let\shiftleft\shiftproofbelowneg
\let\ifwasinsideprooftree\ifinsideprooftree
\insideprooftreetrue
\setbox\proofabove=\hbox\bgroup$\displaystyle 
\let\wereinproofbit\prooftree
\shortenproofleft=0pt \shortenproofright=0pt \proofbelowshift=0pt
\onleftofproofruletrue\penalty1
}

\def\eproofbit{
\ifx    \wereinproofbit\prooftree
\then   \ifcase \lastpenalty
        \then   \shortenproofright=0pt  
        \or     \unpenalty\hfil         
        \or     \unpenalty\unskip       
        \else   \shortenproofright=0pt  
        \fi
\fi
\global\dimen0=\shortenproofleft
\global\dimen1=\shortenproofright
\global\dimen2=\proofrulebreadth
\global\dimen3=\proofbelowshift
\global\dimen4=\proofdotseparation
\global\count255=\proofdotnumber
$\egroup  
\shortenproofleft=\dimen0
\shortenproofright=\dimen1
\proofrulebreadth=\dimen2
\proofbelowshift=\dimen3
\proofdotseparation=\dimen4
\proofdotnumber=\count255
}

\def\proofover{
\eproofbit 
\setbox\proofbelow=\hbox\bgroup 
\let\wereinproofbit\proofover
$\displaystyle
}%
\def\proofoverdbl{
\eproofbit 
\doubleprooftrue
\setbox\proofbelow=\hbox\bgroup 
\let\wereinproofbit\proofoverdbl
$\displaystyle
}%
\def\proofoverdots{
\eproofbit 
\proofdotstrue
\setbox\proofbelow=\hbox\bgroup 
\let\wereinproofbit\proofoverdots
$\displaystyle
}%
\def\proofusing{
\eproofbit 
\setbox\proofrulename=\hbox\bgroup 
\let\wereinproofbit\proofusing
\kern0.3em$
}

\def\endprooftree{
\eproofbit 
  \dimen5 =0pt
\dimen0=\wd\proofabove \advance\dimen0-\shortenproofleft
\advance\dimen0-\shortenproofright
\dimen1=.5\dimen0 \advance\dimen1-.5\wd\proofbelow
\dimen4=\dimen1
\advance\dimen1\proofbelowshift \advance\dimen4-\proofbelowshift
\ifdim  \dimen1<0pt
\then   \advance\shortenproofleft\dimen1
        \advance\dimen0-\dimen1
        \dimen1=0pt
        \ifdim  \shortenproofleft<0pt
        \then   \setbox\proofabove=\hbox{%
                        \kern-\shortenproofleft\unhbox\proofabove}%
                \shortenproofleft=0pt
        \fi
\fi
\ifdim  \dimen4<0pt
\then   \advance\shortenproofright\dimen4
        \advance\dimen0-\dimen4
        \dimen4=0pt
\fi
\ifdim  \shortenproofright<\wd\proofrulename
\then   \shortenproofright=\wd\proofrulename
\fi
\dimen2=\shortenproofleft \advance\dimen2 by\dimen1
\dimen3=\shortenproofright\advance\dimen3 by\dimen4
\ifproofdots
\then
        \dimen6=\shortenproofleft \advance\dimen6 .5\dimen0
        \setbox1=\vbox to\proofdotseparation{\vss\hbox{$\cdot$}\vss}%
        \setbox0=\hbox{%
                \advance\dimen6-.5\wd1
                \kern\dimen6
                $\vcenter to\proofdotnumber\proofdotseparation
                        {\leaders\box1\vfill}$%
                \unhbox\proofrulename}%
\else   \dimen6=\fontdimen22\the\textfont2 
        \dimen7=\dimen6
        \advance\dimen6by.5\proofrulebreadth
        \advance\dimen7by-.5\proofrulebreadth
        \setbox0=\hbox{%
                \kern\shortenproofleft
                \ifdoubleproof
                \then   \hbox to\dimen0{%
                        $\mathsurround0pt\mathord=\mkern-6mu%
                        \cleaders\hbox{$\mkern-2mu=\mkern-2mu$}\hfill
                        \mkern-6mu\mathord=$}%
                \else   \vrule height\dimen6 depth-\dimen7 width\dimen0
                \fi
                \unhbox\proofrulename}%
        \ht0=\dimen6 \dp0=-\dimen7
\fi
\let\doll\relax
\ifwasinsideprooftree
\then   \let\VBOX\vbox
\else   \ifmmode\else$\let\doll=$\fi
        \let\VBOX\vcenter
\fi
\VBOX   {\baselineskip\proofrulebaseline \lineskip.2ex
        \expandafter\lineskiplimit\ifproofdots0ex\else-0.6ex\fi
        \hbox   spread\dimen5   {\hfi\unhbox\proofabove\hfi}%
        \hbox{\box0}%
        \hbox   {\kern\dimen2 \box\proofbelow}}\doll%
\global\dimen2=\dimen2
\global\dimen3=\dimen3
\egroup 
\ifonleftofproofrule
\then   \shortenproofleft=\dimen2
\fi
\shortenproofright=\dimen3
\onleftofproofrulefalse
\ifinsideprooftree
\then   \hskip.5em plus 1fil \penalty2
\fi
}

\begin{document}
\bibliographystyle{plain}
\maketitle

\begin{abstract}
  \zf{} has sentences that quantify over all sets or all ordinals, without restriction.  Some have argued that sentences of this kind lack a determinate meaning.   We propose a set theory called \tops{}, using Natural Deduction, that avoids this problem by speaking only about particular sets.  
\end{abstract}

\section{Introduction}

\zf{} has long been established as a set-theoretic foundation of mathematics, but concerns about its meaningfulness have often been raised.  Specifically, its use of unrestricted quantifiers seems to presuppose an absolute totality of sets.  The purpose of this paper is to present a new set theory called \tops{} that addresses this concern by speaking only about particular sets.  Though weaker than \zf{}, it is adequate for a large amount of mathematics.

To explain the need for a new theory, we begin in Section~\ref{sect:phil} by surveying some basic mathematical beliefs and conceptions.   Section~\ref{sect:lang} presents the language of \tops{}, and Section~\ref{sect:deftops} the theory itself, using Natural Deduction.  Section~\ref{sect:freevar} adapts the theory to allow theorems with free variables.  

Related work is discussed in Section~\ref{sect:related}.  While \tops{} is by no means the first set theory to use restricted quantifiers in order to avoid assuming a totality of sets, it is the first to do so using Natural Deduction, which turns out to be very helpful.  This is apparent when \tops{} is compared to a previously studied system~\cite{Mayberry:globalquant} that proves essentially the same sentences but includes an inference rule that (we argue) cannot be considered truth-preserving.  What saves \tops{} from this pitfall is its use of Natural Deduction.

We sum up the paper's argument in Section~\ref{sect:conclude}.

\section{Motivating \tops{}} \label{sect:phil}
\subsection{Beliefs abut bivalence and definiteness}

Before trying to determine whether a given sentence is true or false, one might wonder whether it has a truth value at all.  A \emph{bivalent} sentence is one that has an objective, determinate truth value (True or False)---regardless of whether anyone can ever know it.  

 Which mathematical sentences are bivalent?  This is a contentious question, and there are various schools of thought.  We shall consider in turn the languages of Peano arithmetic, second and third order arithmetic, and \zf{}.   In each of these, a sentence $\phi$ is built from atomic formulas using connectives and quantifiers.  It is reasonable to say that $\phi$ is bivalent if the range of each quantifier is \emph{definite}, i.e.\ clearly defined, with no ambiguity or haziness.   So we examine each of our languages in the light of this principle.  
\begin{itemize}
\item  Arithmetical sentences, such as the Goldbach Conjecture, quantify over the set $\nats=\setbr{0,1,\ldots}$ of natural numbers.   The conception of $\nats$ is ``most restrictive'', allowing only what is generated by zero and successor.  It is generally considered obvious (though finitists would disagree) that $\nats$ is definite and so arithmetical sentences are bivalent.   
 Let us accept this view and continue.  

\item  Second order arithmetical sentences quantify over $\pset\nats$, or equivalently the set $2^{\nats}$ of $\omega$-sequences of booleans.  While some people (call them \emph{countabilists}) believe in $\nats$ but not in $2^{\nats}$---see e.g.~\cite{Velleman:constructlib,Feferman:predicativism}---the classical view is that, any definite set $X$ has a definite $\omega$-power, i.e.\ set $X^{\nats}$ of $\omega$-sequences in $X$.  Many adherents of this view believe also that $X^{\nats}$ satisfies Dependent Choice, i.e.\ that for any $x \in X$ and entire\footnote{A relation $R$ from $X$ to $Y$ is \emph{entire} when for every $x \smin X$ there is $y \smin Y$ such that $R(x,y)$.} relation
  $R$ from $X$ to $X$, 
 there is $s \smin X^{\nats}$ such that $s_0 = x$ and $\forall n \smin \nats.\, R(s_{n},s_{n+1})$.  These beliefs spring from a ``most liberal'' conception of $\omega$-sequences as consisting of arbitrary choices.   Let us accept them and continue.

\item  Third order arithmetical sentences, such as the Continuum Hypothesis, quantify over $\pset(2^{\nats})$.  While some people (call them \emph{$\omega$-powerists}) believe in all the above but not in $\pset(2^{\nats})$---see e.g.~\cite[Section 3.2.3]{Rathjen:scopefeferman}---the classical view is that any definite set has a definite powerset.  This is equivalent to saying that, for any definite sets $X$ and $Y$, the function set $Y^{X}$ is definite.\footnote{To justify ($\Rightarrow$), take $Y^{X} \eqdef \setbr{r \smin \pset(X \times Y) \mid \forall x \smin X.\, \existsu y \smin Y. \tuple{x,y} \in r}$.  To justify ($\Leftarrow$), take $\pset X \eqdef \setbr{\setbr{x \smin X \mid f(x) =1} \mid f \smin \setbr{0,1}^{X}}$.}  Many adherents of this view believe also that $Y^{X}$ satisfies the Axiom of Choice, i.e.\ that for any entire relation $R$ from $X$ to $Y$,
 there is $f \smin Y^{X}$ such that $\forall x \smin X.\, R(x , f(x))$.   These beliefs (to paraphrase Bernays~\cite{Bernays:platonism}) spring from a ``most liberal'' conception of functions or subsets as consisting of arbitrary choices.  Let us accept them and continue.  

\item Recall next the definition of the \emph{cumulative hierarchy}.  It associates to each ordinal $\alpha$ a set $V_{\alpha}$ by transfinite recursion:
  \begin{spaceout}{lrcl}
&  V_{\alpha+1} & \eqdef & \pset V_{\alpha} \\
\text{If $\alpha$ is a limit, } & V_{\alpha } & \eqdef & \bigcup_{\beta < \alpha} V_{\beta}  \\
\text{In particular, } & V_{0} & \eqdef & \emptyset 
\end{spaceout}%
For example, consider the ordinal $\omega_1$.  (It  can be implemented as the set of well-ordered subsets of $\nats$, modulo isomorphism.)  Is $V_{\omega_1}$ definite? The classical view---that, for any definite ordinal $\alpha$, the set $V_{\alpha}$ is definite---may be argued by induction on $\alpha$ as follows.   For the successor case: any definite set has a definite powerset.  For the limit case: given a definite set $I$, and for each $i \in I$ a definite set $A_i$, surely the union $\bigcup_{i \in I} A_i$ is definite.   Let us accept this argument and continue. 
\item  \zf{} has sentences that quantify over the ordinals.   To emphasize: over \emph{all} the ordinals, however large.  That is a ``most liberal'' conception.
Is it definite?  Arguably not, because of the Burali-Forti paradox: the totality of all ordinals is itself well-ordered, with an order-type too large to belong to it.
This powerful argument led Parsons to describe the language of \zf{} as ``systematically ambiguous''~\cite{Parsons:setsandclasses}, and was taken up by Dummett~\cite[pages 316--317]{Dummett:fregephilmath}:
  \begin{quotation}
    The Burali-Forti paradox ensures
that no definite totality comprises everything intuitively recognisable
as an ordinal number,  where a definite totality is one over which quantification always yields a statement that is determinately true
or false.
\end{quotation}
However, it was rejected by Boolos~\cite{Boolos:whencecontra} and Cartwright~\cite{Cartwright:speakingevery}, who held that quantifying over the ordinals does not imply that they form a totality.\footnote{Even for those who hold this view, it is hard to deny that quantifying over \emph{classes} of ordinals, as in Morse-Kelley class theory, implies that the ordinals form a totality.  And harder still for higher levels of class.}
\end{itemize}

\noindent To summarize: we have accepted the classical view that (for example) the set $V_{\omega_1}$ is definite and, by the Axiom of Choice, well-orderable.  But we have not accepted that quantification over the ordinals is definite.  Since we believe in many \emph{particular} sets, such as $V_{\omega_{1}}$, our view may be called \emph{particularism}.  Other names given to it, or to similar views, are ``restricted platonism''~\cite{Bernays:platonism,McNaughton:conceptualschemes},  ``liberal intuitionism''~\cite{Pozsgay:liberal} and ``Zermelian potentialism''~\cite{AntosFriedmanHonzikTernullo:hyperun}
---see the discussion in~\cite[Section 4.3]{Friedman:evidencehyperun}.

\subsection{\zf{} and purity}

Recall next that \zf{} adopts the following assumptions.
  \begin{itemize} 
  \item \emph{Everything is a set.}  This rules out, for example, urelements and primitive ordered pairs.
  \item \emph{Membership is well-founded.}  This rules out, for example, sets $a$ such that $a = \setbr{a}$, known as ``Quine atoms''.
  \end{itemize}
We call them \emph{purity assumptions}.  They can be combined into the statement \emph{Everything is a well-founded pure set}, where ``pure'' means that every element, and every element of an element, etc., is a set.  The significance of these assumptions depends on how we read quantifiers in \zf{}.
  \begin{itemize}
  \item The \emph{full-blown} interpretation of $\forall x$ is ``for any thing $x$ whatsoever''.
     On this reading, \zf{} denies the existence of urelements and Quine atoms, which makes it unsound if an urelement or Quine atom does in fact exist (whatever that means).
  \item The \emph{pure} interpretation of $\forall x$ is ``for any well-founded pure set $x$''.  On this reading, \zf{} neither denies nor affirms the existence of urelements and Quine atoms.  It simply refrains from speaking about them.
  \end{itemize}
 Each interpretation raises a question.
 \begin{itemize}
    \item Is quantification over all things definite?  That seems implausible.   For example, it is hard to believe that the question ``Are there precisely seven urelements?'' has an objectively correct answer.
    \item Is quantification over the well-founded pure sets definite?  The answer is yes if and only if  quantification over the ordinals is definite.\footnote{To justify ($\Leftarrow$), note that 
  the quantifier ``for any well-founded pure set $x$'' is equivalent to ``for any ordinal $\alpha$ and any $x \smin V_{\alpha}$'', since any well-founded pure set belongs to $V_{\alpha}$, where $\alpha$ is the successor of its rank.  To justify ($\Rightarrow$), use von Neumann's  implementation of the ordinals as the transitive pure sets that are well-ordered by membership.}
  \end{itemize}
A \emph{full-blown totalist} believes in quantification over all things.  A \emph{pure totalist} does not believe this, but does believe in quantification over the well-founded pure sets, or equivalently over the ordinals.






\subsection{Set theory for particularists}


We can now state the problem.  \zf{} is a \emph{set theory for pure totalists}, in the sense that pure totalists (reading it purely) would view its sentences as bivalent and its theorems as true.   What would be a set theory for particularists (i.e.\ us), in this sense?  Such a theory would enable us to state and prove facts about particular sets, e.g.\ ``The set $V_{\omega_1}$ is well-orderable'', but prevent us from forming sentences that quantify over all sets ({``Every set is well-orderable''}) or all ordinals ({``Ordinal addition is associative''}) or all things (``Everything is equal to itself'').



This paper presents a system meeting these requirements called \emph{TOPS}---short for ``Theory Of Particular Sets''.  It differs from ZFC in two substantive ways: it avoids unrestricted quantification, and does not adopt the purity assumptions.    The two points of difference are, in principle, orthogonal, but the purity assumptions are needed in ZFC to make the quantifiers intelligible to pure totalists.  In a system where quantifiers are always restricted, they are not needed.

\section{The language of \tops{}} \label{sect:lang}
\subsection{Informal account} \label{sect:notation}

We begin by informally introducing the \tops{} notation.  To this end, let $\totall$ be the totality of all things.  It may contain urelements, Quine atoms, inaccessible cardinals, measurable cardinals, etc.  Let $\allsets$ (a subcollection of $\totall$) be the totality of all sets. 

The \tops{} notation is listed as follows; most of it is familiar.
\begin{enumerate}
\item Let $a$ be a thing.  We write $\isset{a}$ to say that $a$ is a set.
\item Let $A$ be a set and $b$ a thing.  We write $b \in A$ to say that $b$ is an element of $A$.
\item Let $A$ and $B$ be sets.  We say that $A$ is \emph{included in} $B$, written $A \subseteq B$, when every element of $A$ is in $B$.
   \item We write $\bbot$ for the empty set.
\item Let $A$ and $B$ be sets.  We write $A \cup B$ for the union of $A$ and $B$.
\item Let $A$ be a set and $(B(a))_{a \in A}$ a family of sets.  We write $\bigcup_{a \in A} B(a)$ for the union of the family.
\item Let $A$ be a set and $P$ a predicate on it.  We write $\setbr{x \smin A \mid P(x)}$ for the set of $x \in A$ satisfying $P$.  
\item Let $a$ be a thing.
  \begin{itemize}
  \item The \emph{singleton} $\setbr{a}$ is the set whose sole element is $a$.
  \item For an assertion $\phi$, the \emph{conditional singleton} $\condsing{a}{\phi}$ is $\setbr{a}$ if $\phi$ is true and $\emptyset$ otherwise.
\end{itemize}
\item Let $A$ be a set.
  \begin{itemize}
  \item The range of a function $f \colon A \to \totall$ is written
    $\setbr{f(x) \mid x \smin A}$.
   \item For a predicate 
    $P$ on $A$, the range of a function
    $f \colon \setbr{x \in A \mid P(x)} \to \totall$ is written
    $\setbr{f(x) \mid x \smin A \mid P(x)}$.
\end{itemize}
\item Let $A$ be a singleton set.  Its unique element is written $\uniqel{A}$.
\item Let $a$ be a thing and $F \colon \totall \to \totall$ a function, which may be written $x.F(x)$.  The \emph{iterative reach} of $F$ on $a$, written $\itset{a}{F}$, is the set $\setbr{f^{n}(a) \mid n \in \nats}$.
\item  Let $A$ be a set and $R$ a relation on $A$, which may be written $x,y.R(x,y)$. 
  \begin{itemize}
  \item We write $\wfint{a}{A}{R}$ when $a$ is an $R$-well-founded element of $A$.
This property is generated inductively by the following rule: for $a \in A$, if every $b \in A$ such that $R(b,a)$ is $R$-well-founded, then so is $a$.  
By Dependent Choice, an element $a$ of $A$ is $R$-well-founded iff there is no sequence $(a_n)_{n \in \nats}$ of elements of $A$ such that $a_0 = a$ and $\forall n \smin \nats.\,R(a_{n+1},a_n)$.
\item Let $F \colon A \times \allsets \to \totall$ be a function, which may be written $z,Y.F(z,Y)$.  For any $R$-well-founded element $a$ of $A$, we use \emph{well-founded recursion} to define a thing written $\prl{A}{R}{F}{a}$.  Here is the recursive definition:
  \begin{eqnarray*}
        \prl{A}{R}{F}{a} & \eqdef & F(a,\setbr{\prl{A}{R}{F}{x} \mid x \smin A \mid R(x,a)}) 
   \end{eqnarray*}
\end{itemize}
\item Let $A$ be a set.  Its \emph{powerset}, written $\pset A$, is the set of all subsets of $A$.
\end{enumerate}

To ensure that our notation is defined in all cases, we adopt the following conventions.
\begin{enumerate}
\item \emph{Nonsense denotes the empty set
.} 
For a set $A$, we take
\begin{itemize}
\item $\uniqel{A}$ to be $\bbot$ if $A$ is not singleton
\item 
  $\prl{A}{R}{F}{a}$ to be $\bbot$ if $a$ is not an
  $R$-well-founded element of $A$.
\end{itemize}
For an assertion $\phi$ and thing $a$, we say that $\phi$ is a \emph{Nonsense Convention prerequisite} for $a$, written $\nonsenseconv{\phi}{a}$, when $\neg \phi$ implies $a = \emptyset$.   This can also be stated in a positive way: either $\phi$ is true, or $a$ is a set such that, if it is inhabited, then $\phi$ is true. 



\item \emph{Non-sets are treated like the empty set.}  Wherever the notation expects a set, any non-set provided is tacitly replaced by $\emptyset$.   For example, if $a$ is not a set, we take
  \begin{itemize}
  \item the statement $b \in a$ to be false
  \item the statement $a \subseteq b$ to be true 
   \item $\pset a$ to be $\setbr{\bbot}$.
\end{itemize}
Likewise, for a set $A$ and family of things $(b(a))_{a \in A}$, we take $\bigcup_{a \in A} b(a)$ to be $\bigcup_{a \in C} b(a)$, where $C$ is the set of $a \in A$ such that $b(a)$ is a set.



\end{enumerate}

\subsection{Formal syntax}

We now present the syntax formally.  Let $\vars$ be an infinite set of variables, written $x,y,X,Y,\ldots$.
The syntax of terms and that of propositional formulas are mutually inductively defined in Figure~\ref{fig:przfcexp}.   We usually use a lowercase letter to suggest a thing, an uppercase one to suggest a set, and a calligraphic one (e.g.\ $\cata$) to suggest a set of sets.  
But this is just for the sake of readability; the system does not distinguish these.

\begin{figure}
  \centering
   \begin{spaceout}{lrcl}
\text{Terms} & r,s,t,A,B & \bnfgo & x \bnf \bbot \bnf A \cup B \bnf \bigcup_{a \in A} B(a)  \bnf \condsing{r}{\phi} \bnf  \uniqel{A}
\\
& & & \bnf \itset{s}{x.r}  \bnf \prl{A}{x,y.\phi}{z,Y.r}{s}  \bnf \pset A \\
\multicolumn{2}{l}{\text{Propositional formulas}  \quad\quad \phi,\psi} & \bnfgo & \False \bnf \True \bnf \phi \vee \psi   \bnf \phi \wedge \psi  \bnf \phi \Rightarrow \psi \bnf \neg \phi \\
& & & \bnf \exists x \smin A.\,\phi 
\bnf \forall x \smin A.\,\phi \bnf  s = t  \\
& & & \bnf \isset{A}  \bnf s \in A 
\bnf \wfint{s}{A}{x,y.\phi}
  \end{spaceout}
  \caption{Syntax of \tops{}}
  \label{fig:przfcexp}
\end{figure}

A \emph{declaration context} $\gamma$ is a finite subset of $\vars$.  We define $\mathsf{nil} \eqdef \emptyset$ (it may be written as nothing) and $\gamma,x \eqdef \gamma \cup \setbr{x}$.   We write
\begin{itemize}
\item $\sht{\gamma}{t}$ to say that $t$ is a term over $\gamma$, i.e.\ a term whose free variables are all in $\gamma$
\item $\shp{\gamma}{\phi}$  to say that $\phi$ is a propositional formula over $\gamma$.
\end{itemize}
These two judgements are defined by mutual  induction in Figure~\ref{fig:termform}.  

\begin{figure}
  \begin{displaymath}
    \begin{array}{ccccccc}
      \begin{prooftree}
        \strut
        \using (x \in \gamma)
        \justifies
        \sht{\gamma}{x}
      \end{prooftree}
        & & 
\begin{prooftree}
                        \strut
                        \justifies
                        \sht{\gamma}{\bbot}
                      \end{prooftree}
 & &
      \begin{prooftree}
        \sht{\gamma}{A} \betwixt \sht{\gamma}{B}
        \justifies
        \sht{\gamma}{A \cup B}
      \end{prooftree} & & 
                          \begin{prooftree}
                            \sht{\gamma}{A} \betwixt \sht{\gamma,x}{B}
                            \justifies
                            \sht{\gamma}{\bigcup_{x \in A}B}
                          \end{prooftree} \\ \\
      \begin{prooftree}
        \sht{\gamma}{r} \betwixt \shp{\gamma}{\phi}
        \justifies
        \sht{\gamma}{\condsing{r}{\phi}}
      \end{prooftree} & & 
\begin{prooftree}
     \sht{\gamma}{A}
     \justifies                                                             \sht{\gamma}{\uniqel{A}}
                        \end{prooftree} & & 
  \begin{prooftree} 
\sht{\gamma}{s} \betwixt \sht{\gamma,x}{r}
        \justifies
        \sht{\gamma}{\itset{s}{x.r}}
      \end{prooftree} \\ \\
  \multicolumn{3}{c}{    \begin{prooftree}
                          \sht{\gamma}{A} \betwixt \shp{\gamma,x,y}{\phi} \betwixt \sht{\gamma,z,Y}{r} \betwixt \sht{\gamma}{s}
                          \justifies
                          \sht{\gamma}{\prl{A}{x,y.\phi}{z,Y.r}{s}}
                        \end{prooftree} }
                  & &
      \begin{prooftree}
        \sht{\gamma}{A}
        \justifies
        \sht{\gamma}{\pset A}
      \end{prooftree} \\ \\
  \begin{prooftree}
                          \strut
                          \justifies
                          \shp{\gamma}{\False}
                        \end{prooftree} & & 
\begin{prooftree}
                          \strut
                          \justifies
                          \shp{\gamma}{\True}
                        \end{prooftree} & &
 \begin{prooftree}
        \shp{\gamma}{\phi} \betwixt \shp{\gamma}{\psi}
        \justifies
        \shp{\gamma}{\phi \vee \psi}
      \end{prooftree} & & 
 \begin{prooftree}
        \shp{\gamma}{\phi} \betwixt \shp{\gamma}{\psi}
        \justifies
        \shp{\gamma}{\phi \wedge \psi}
      \end{prooftree} \\ \\
 \begin{prooftree}
        \shp{\gamma}{\phi} \betwixt \shp{\gamma}{\psi}
        \justifies
        \shp{\gamma}{\phi \Rightarrow \psi}
      \end{prooftree} & & 
 \begin{prooftree}
        \shp{\gamma}{\phi}
        \justifies
        \shp{\gamma}{\neg \phi}
      \end{prooftree} & &  
 \begin{prooftree}
                          \sht{\gamma}{A} \betwixt \shp{\gamma,x}{\phi}
                          \justifies
                          \shp{\gamma}{\exists x \smin A.\,\phi}
                         \end{prooftree} & & 
\begin{prooftree}
                          \sht{\gamma}{A} \betwixt \shp{\gamma,x}{\phi}
                          \justifies
                          \shp{\gamma}{\forall x \smin A.\,\phi}
                        \end{prooftree} \\ \\
\begin{prooftree}
        \sht{\gamma}{s} \betwixt \sht{\gamma}{t}
        \justifies
        \shp{\gamma}{s = t}
      \end{prooftree} 
& & 
      \begin{prooftree}
        \sht{\gamma}{A}
        \justifies
        \shp{\gamma}{\isset{A}}
      \end{prooftree} & &
\begin{prooftree}
        \sht{\gamma}{s} \betwixt \sht{\gamma}{A}
        \justifies
        \shp{\gamma}{s \in  A}
      \end{prooftree} & &
\begin{prooftree}
        \sht{\gamma}{s} \betwixt \sht{\gamma}{A} \betwixt \shp{\gamma,x,y}{\phi}
        \justifies
        \shp{\gamma}{\wfint{s}{A}{x,y.\phi}}
      \end{prooftree}
    \end{array}
  \end{displaymath}
  \caption{Terms and propositional formulas in context}
  \label{fig:termform}
\end{figure}


As usual, we identify terms up to renaming of bound variables ($\alpha$-equivalence).   We write $\phi[t/x]$ and $s[t/x]$ for the capture-avoiding substitution of $t$ for $x$ in $\phi$ and in $s$ respectively.  

A \emph{sentence} is a closed propositional formula, so $\shp{}{\phi}$ says that $\phi$ is a sentence.

\section{Definition of \tops{}} \label{sect:deftops}
\subsection{Sequents} \label{sect:sequent}

We formulate \tops{} using Natural Deduction~\cite{PelletierHazen:histnatded}, which is a convenient framework for proofs in general, but especially for proofs involving restricted quantification and nested dependencies, as dependent type theory has shown~\cite{MartinLoef:inttypetheory}.  

In the middle of a Natural Deduction proof, an assertion is made subject to some hypotheses and variable declarations.  This information can be presented as a \emph{sequent}.  For example, given formulas and terms
\begin{spaceout}{lcl}
 & \vdashp & \phi_0 \\
 & \vdasht & A \\ 
x & \vdashp & \phi_1 \\
x & \vdasht & B \\
x,y & \vdasht & C \\
x,y,z & \vdashp & \psi 
\end{spaceout}%
the following is a sequent:
\begin{displaymath} 
  \phi_0 \lco x \colin A \lco \phi_1 \lco y \colin B \lco z \colin C\ \vdash\ \psi
\end{displaymath}
It is read: ``Assuming $\phi_0$, and $x$ in $A$, and $\phi_1$, and $y$ in $B$, and $z$ in $C$, we assert $\psi$.''  Thus it has the same meaning as the sentence
\begin{displaymath}
  \phi_0 \Rightarrow \forall x \smin A.\, (\phi_1 \Rightarrow \forall y \smin B.\,\forall z \smin C.\,\psi)
\end{displaymath}
The part of a sequent that appears to the left of the $\vdash$ symbol is called a \emph{logical context}.  Let us now make this precise.


\begin{definition}
We define the set of \emph{logical contexts}, and to each of these we associate
\begin{itemize}
\item a declaration context $\dvar{\Gamma}$
\item and a list $\hyps{\Gamma}$ of formulas over $\dvarg{\Gamma}$, called \emph{hypotheses},
\end{itemize}
inductively as follows.
\begin{itemize}
\item $\mathsf{nil}$ is a logical context (it may be written as nothing), with
\begin{eqnarray*}
    \dvar{\mathsf{nil}} & \eqdef & \mathsf{nil}                 \\
    \hyps{\mathsf{nil}} & \eqdef & \text{the empty list}
  \end{eqnarray*}
\item If $\Gamma$ is a logical context and $x \not\in \dvarg{\Gamma}$ and $\sht{\dvarg{\Gamma}}{A}$, then $\Gamma\lco x\colin A$ is a logical context, with
 \begin{eqnarray*}
    \dvar{\Gamma\lco  x \colin A} & \eqdef & \dvar{\Gamma}\lco x \\
    \hyps{\Gamma\lco x \colin A} & \eqdef & \hyps{\Gamma}\lco x \smin A
  \end{eqnarray*}
\item  If $\Gamma$ is a logical context and $\shp{\dvar{\Gamma}}{\phi}$, then $\Gamma\llco \phi$ is a logical context, with
  \begin{eqnarray*}
    \dvar{\Gamma\llco \phi} & \eqdef & \dvar{\Gamma} \\
    \hyps{\Gamma\llco \phi} & \eqdef & \hyps{\Gamma}\llco\phi
  \end{eqnarray*}
\end{itemize}
\end{definition}

\begin{definition}
  A \emph{sequent}, written $\Gamma \vdash \psi$, consists of
  \begin{itemize}
  \item a logical context $\Gamma$
  \item and a formula $\psi$ over $\dvar{\Gamma}$, called the \emph{conclusion}.
  \end{itemize}
\end{definition}
Below (Definition~\ref{def:formtoseq}) we represent a sentence $\phi$ as the sequent $ \vdash \phi$.

\subsection{Abbreviations}

We shall introduce some abbreviations to aid readability.  

For $n \in\nats$, an \emph{$n$-ary abstracted term} is a term $F$ that may use additional variables $z_0,\ldots,z_{n-1}$ not in $\vars$. It is said to be over $\gamma$ when all its free variables other than these are in $\gamma$.   Given terms $\vec{t} = t_0,\ldots,t_{n-1}$ we write 
\begin{eqnarray*}
  F(\vec{t}) & \eqdef & F[t_0/z_0,\ldots,t_{n-1}/z_{n-1}]
\end{eqnarray*}
If $F$ and $\vec{t}$ are over $\gamma$ then so is $F(\vec{t})$.  We likewise define $P(\vec{t})$, where $P$ is an \emph{$n$-ary abstracted propositional formula}. We often abbreviate
$\vec{x}.\,F(\vec{x})$ as $F$ and $\vec{x}.\,P(\vec{x})$ as $P$.



We write $A_0 \cup \cdots \cup A_{n-1}$ and $\phi_{0} \vee \cdots \vee \phi_{n-1}$ and $\phi_{0} \wedge \cdots \wedge \phi_{n-1}$ in the usual way.  For $n=0$ these are $\emptyset$ and $\False$ and $\True$ respectively, and for $n \geqslant 3$ we choose some arrangement of parentheses.

We make the following abbreviations.   
\begin{spaceout}{rcl}
 \phi \Leftrightarrow \psi & \eqdef & (\phi \Rightarrow \psi) \wedge (\psi \Rightarrow \phi) 
 \\
 \existsu x \smin A.\, P(x) & \eqdef & \exists x \smin A.\, (P(x) \wedge \forall y \smin A.\, (P(y) \Rightarrow y=x))  \betwixt \betwixt  
 \\
 A \subseteq B & \eqdef & \forall x \smin A.\,x \in B \\
 \setbr{x \in A \mid P(x)} & \eqdef & \bigcup_{x \in A} \condsing{x}{P(x)} \\
 \setbr{r} & \eqdef & \condsing{r}{\True} \\
 \setbr{r_0,\ldots,r_{n-1}} & \eqdef & \setbr{r_0} \cup \cdots \cup \setbr{r_{n-1}} \\
  \setbr{F(x) \mid x \smin A} & \eqdef & \bigcup_{x \in A} \setbr{F(x)} \\
 \setbr{F(x,y) \mid x \smin A, y \smin B(x)} & \eqdef & \bigcup_{x \in A}\bigcup_{y \in B(x)} \setbr{F(x,y)} \betwixt \text{etc.} \\
 \setbr{F(x) \mid x \smin A \mid P(x)} & \eqdef & \bigcup_{x \in A} \condsing{F(x)}{P(x)} \\
 \setbr{F(x,y) \mid x \smin A, y \smin B(x) \mid P(x,y)} & \eqdef & \bigcup_{x \in A}\bigcup_{y \in B(x)} \condsing{F(x,y)}{P(x,y)} \betwixt \text{etc.} \\
 A \cap B & \eqdef & \setbr{x \smin A \mid x \in B} \\
  \rotiota x \smin A.\, P(x) & \eqdef & \uniqel{\setbr{x \smin A \mid P(x)}}  \\
\bigcup \cata & \eqdef & \bigcup_{X \in \cata} X  \\
\nonsenseconv{\phi}{r} & \eqdef & \phi \; \vee\; (\isset{r} \wedge \forall x \in r.\, \phi)
\end{spaceout}%
Here and throughout the paper, fresh variables are used for binding.  For example, in the definition of $\existsu x \smin A.\, P(x)$ above, $y$ must be fresh for $P$ and not be $x$.

\subsection{Provability} \label{sect:proveseq}


\tops{} consists of logical rules (Figure~\ref{fig:logrules}) and axiom schemes (Figure~\ref{fig:axiomschemes}), which together define provability of sequents.  Each rule and scheme refers to a logical context $\Gamma$.  All (abstracted) terms and formulas mentioned are assumed to be over $\dvar{\Gamma}$.  The logical rules are the usual  introduction ($\Intro$) and elimination ($\Elim$) rules of Natural Deduction, adapted to restricted quantification.    As there is no standard set of Natural Deduction rules for negation, the $\neg$Quodlibet and $\neg$Cases rules have been chosen for convenience and symmetry.

\begin{figure}
  \begin{displaymath}
    \begin{array}{ccccc}
      \begin{prooftree}
\phi \text{ appears in } \hyps{\Gamma}
        \using \parenth{\text{Hypothesis}}
        \justifies
        \Gamma \vdash \phi
      \end{prooftree} & & 
      \begin{prooftree}
        \Gamma \vdash \False\using \parenth{\False\,\Elim} \justifies \Gamma \vdash \psi
      \end{prooftree} & &
\begin{prooftree}
                          \strut \using \parenth{\True\,\Intro} \justifies \Gamma \vdash \True
                        \end{prooftree} \\ \\
                        \begin{prooftree}
                          \Gamma \vdash \phi \using \parenth{\vee\Intro L} \justifies
                          \Gamma \vdash \phi\vee \phi'
                        \end{prooftree} & &
                                          \begin{prooftree}
                                            \Gamma \vdash \phi'\using \parenth{\vee\Intro R}
                                            \justifies \Gamma
                                            \vdash \phi\vee \phi'
                                          \end{prooftree} & 
 \multicolumn{2}{l}{    \begin{prooftree}
        \Gamma \vdash \phi \vee \phi' \betwixt \Gamma\llco\phi
        \vdash \psi \betwixt \Gamma\llco\phi' \vdash \psi \using \parenth{\vee\Elim} \justifies
        \Gamma \vdash \psi
      \end{prooftree} } 
                         \\ \\
                                          \begin{prooftree}
                                            \Gamma \vdash \phi
                                            \betwixt \Gamma
                                            \vdash \phi' \using \parenth{\wedge\Intro} \justifies
                                            \Gamma \vdash
                                            \phi\wedge \phi'
                                          \end{prooftree} & &
      \begin{prooftree}
        \Gamma \vdash\phi \wedge \phi' \using \parenth{\wedge\Elim L} \justifies \Gamma
        \vdash \phi
      \end{prooftree} & &
                        \begin{prooftree}
                          \Gamma \vdash\phi \wedge \phi'
                          \using \parenth{\wedge\Elim R} \justifies \Gamma \vdash \phi'
                        \end{prooftree} \\ \\
                                          \begin{prooftree}
                                            \Gamma\llco\phi
                                            \vdash\psi \using \parenth{\Rightarrow\!\Intro} \justifies
                                            \Gamma \vdash \phi
                                            \Rightarrow \psi
                                          \end{prooftree} & &
      \begin{prooftree}
        \Gamma \vdash \phi \Rightarrow \psi \betwixt
        \Gamma \vdash \phi \using \parenth{\Rightarrow\!\Elim} \justifies \Gamma \vdash \psi
      \end{prooftree} \\ \\
                        \begin{prooftree}
                          \Gamma \vdash \phi \betwixt
                          \Gamma \vdash \neg \phi \using \using \parenth{\neg\text{Quodlibet}} \justifies
                          \Gamma \vdash \psi
                        \end{prooftree} & &
                                          \begin{prooftree}
                                            \Gamma\llco\phi \vdash
                                            \psi \betwixt
                                            \Gamma\llco\neg \phi
                                            \vdash \psi \using \parenth{\neg\text{Cases}} \justifies
                                            \Gamma \vdash \psi
                                          \end{prooftree} \\ \\
      \begin{prooftree}
        \Gamma \vdash t \in A \betwixt \Gamma \vdash
        P(t) \using \parenth{\exists_{\in}\Intro} \justifies \Gamma \vdash \exists x\smin A.\,
        P(x)
      \end{prooftree} & &
                       \multicolumn{2}{l}{ \begin{prooftree}
                          \Gamma \vdash \exists x \smin A.\, P(x)
                          \betwixt \Gamma\llco x \colin A \llco P(x) \vdash
                          \psi \using \parenth{\exists_{\in}\Elim} \justifies
                          \Gamma \vdash \psi
                        \end{prooftree} } \\ \\
      \begin{prooftree}
        \Gamma\llco x \colin A \vdash P(x) 
        \using \parenth{\forall_{\in}\Intro} \justifies \Gamma \vdash \forall x \smin A.\, P(x)
      \end{prooftree} & &
                        \begin{prooftree}
                          \Gamma \vdash \forall x \smin A.\, P(x)
                          \betwixt \Gamma \vdash t \in A
                          \using \parenth{\forall_{\in}\Elim} \justifies \Gamma \vdash P(t)
                        \end{prooftree} 
 \\ \\
\begin{prooftree}
                             \strut
                             \using \parenth{=\!\Intro} \justifies
                             \Gamma \vdash s=s
                           \end{prooftree} & & 
                                               \begin{prooftree}
                                                 \Gamma \vdash r=s \betwixt \Gamma \vdash P(r)
                                                \using \parenth{=\!\Elim} \justifies
                                                 \Gamma \vdash P(s)
                                               \end{prooftree} 
    \end{array}
  \end{displaymath}
  \caption{Logical rules}
  \label{fig:logrules}
\end{figure}

\begin{figure}
  \begin{spaceout}{llll}
   \multicolumn{4}{l}{\textb{Sethood}} \\
 \textp{Set Membership} & \Gamma & \vdash & s \in A \ \Rightarrow\  \isset{A} \\
 \textp{Set Extensionality} & \Gamma & \vdash & (\isset{A} \wedge  \isset{B} \wedge A \subseteq B \wedge B \subseteq A)\ \Rightarrow\ A=B \interline
 \multicolumn{4}{l}{\textb{Empty set}} \\
 \textp{$\bbot$\bd Set} & \Gamma & \vdash & \isset{\bbot} \\
 \textp{$\bbot$\bd Element} &  \Gamma & \vdash &  t \in \bbot \ifflong \False \interline
 \multicolumn{4}{l}{\textb{Binary union}} \\
\textp{$\cup$\bd Set} &   \Gamma & \vdash & \isset{A \cup B} \\
\textp{$\cup$\bd Element} &   \Gamma & \vdash &  t \in A \cup B \ifflong (t \in A \vee t \in B) \interline
\multicolumn{4}{l}{\textb{Indexed union}} \\
\textp{$\bigcup_{\in}$\bd Set} &  \Gamma & \vdash & \isset{\bigcup_{x \in A} B(x)} \\
\textp{$\bigcup_{\in}$\bd Element} & \Gamma & \vdash &  t \in \bigcup_{x \in A} B(x) \ifflong \exists x \smin A.\, t \in B(x) \interline
\multicolumn{4}{l}{\textb{Conditional singleton}} \\
\textp{CS\bd Set} &   \Gamma & \vdash & \isset{\condsing{r}{\phi}} \\
\textp{CS\bd Element} &  \Gamma & \vdash &  t \in \condsing{r}{\phi} \ifflong (t=r \wedge \phi) \interline
\multicolumn{4}{l}{\textb{Unique element}} \\
\textp{$\mathsf{UE}$\bd $\mathsf{NC}$} &  \Gamma & \vdash & \nonsenseconv{(\existsu x \smin A.\, \True)}{\uniqel{A}} \\
\textp{$\mathsf{UE}$\bd Specification} & \Gamma & \vdash & (\existsu x \smin A.\, \True)\  \Rightarrow\  \uniqel{A} \in A 
\interline
\multicolumn{4}{l}{\textb{Iterative reach}} \\
  \textp{$\mathsf{IR}$\bd Set} &  \Gamma & \vdash & \isset{\itset{s}{F}} \\
 \textp{$\mathsf{IR}$\bd Base Generation} & \Gamma & \vdash & s \in \itset{s}{F} \\
  \textp{$\mathsf{IR}$\bd Step Generation} & \Gamma & \vdash & t \in \itset{s}{F}\ \Rightarrow\ F(t)\ \in \itset{s}{F} \\
  \textp{$\mathsf{IR}$\bd Induction} & \Gamma & \vdash &  (P(s)\, \wedge\, (\forall x \smin \itset{s}{F}.\,(P(x)  \Rightarrow P(F(x))))\, \wedge\, t \in \itset{s}{F})\ \Rightarrow\ P(t) \interline
 \multicolumn{4}{l}{\textb{Well-founded elementhood}} \\
\textp{$\mathsf{W}_{\in}$\bd Generation} &   \Gamma & \vdash & (r \in A \wedge \forall y \smin A.\,R(y,x) \Rightarrow \wfint{y}{A}{R}) \ \Rightarrow \  \wfint{r}{A}{R} \\
\textp{$\mathsf{W}_{\in}$\bd Induction} & \Gamma & \vdash & ((\forall x \smin A.\, (\forall y \smin A.\,(R(y,x) \Rightarrow P(y))) \Rightarrow P(x))  \, \wedge \, \wfint{t}{A}{R})\ \Rightarrow\ P(t) \interline
\multicolumn{4}{l}{\textb{Well-founded recursion}} \\
\textp{$\mathsf{WR}_{\in}$\bd $\mathsf{NC}$} &   \Gamma & \vdash & \nonsenseconv{\wfint{s}{A}{R}}{\prl{A}{R}{F}{s}}  \\
\textp{$\mathsf{WR}_{\in}$\bd Specification} &  \Gamma & \vdash &  \wfint{s}{A}{R} \ \Rightarrow\  \prl{A}{R}{F}{s}  \,  = \,   F(s,\setbr{\prl{A}{R}{F}{x} \mid x \smin A \mid R(x,s)}) 
\interline
\multicolumn{4}{l}{\textb{Powerset}} \\
\textp{$\pset$\bd Set} & \Gamma & \vdash & \isset{\pset A} \\
\textp{$\pset$\bd Element} &  \Gamma & \vdash & 
B \in \pset A \ifflong (\isset{B} \wedge B \subseteq A)
\interline
\multicolumn{4}{l}{\textb{Choice}
}  \\
\textp{Choice} &   \Gamma & \vdash & ((\forall x \smin A.\, \forall w \smin A. \forall y \smin B(x) \cap B(w).\, x=w)\, \wedge\, (\forall x \smin A.\,\exists y \smin B(x).\,\True)) \\
& & & \quad \Rightarrow \ \exists Y \smin \pset \bigcup_{x \in A} B(x).\,\forall x \smin A.\, \existsu y \smin B(x).\,y \in Y
\end{spaceout}%
\caption{Axiom schemes of TOPS}
\label{fig:axiomschemes}
\end{figure}

To complete the definition of Closed \tops{}, we define provability of sentences.
\begin{definition} \label{def:formtoseq}
  A \emph{theorem} is a sentence $\phi$ such that the sequent $ \vdash \phi$ is provable.
\end{definition}




\subsection{Subsystems} \label{sect:subsyst}

Certain subsystems of \tops{} may be of special interest.  For example:
\begin{itemize}
\item the \emph{arithmetical} fragment, which excludes  $\mathsf{W}$, $\mathsf{WR}$, powerset and Choice
\item the \emph{$W$-arithmetical} fragment, which excludes powerset and Choice
\item the \emph{intuitionistic} fragment, which excludes negation and Choice
  \item \tops{} without Choice.
  \end{itemize}
Our formulation of \tops{} has been arranged so as to keep these subsystems self-contained.  
However, this is a matter of taste;  other presentations are possible and may have their own advantages.

\section{Open \tops{}} \label{sect:freevar}

As we have seen, \tops{} allows us to prove \emph{sentences}.  For a declaration context $\gamma$, we shall give a variant of \tops{} called \emph{\tops{} over $\gamma$} allowing us to prove formulas over $\gamma$. These systems are collectively called \emph{Open \tops{}}.

We first define \emph{logical context over $\gamma$} and \emph{sequent over $\gamma$} the same way as in Section~\ref{sect:sequent}, except that we replace $\mathsf{Decl}$ with $\mathsf{Decl}_{\gamma}$.  The definition of the latter differs only in the following clause:
\begin{eqnarray*}
  \mathsf{Decl}_{\gamma}({\mathsf{nil}}) & \eqdef & \gamma       
\end{eqnarray*}
Provability of sequents is defined as in Section~\ref{sect:proveseq}.  Finally, a \emph{theorem over $\gamma$} is a formula $\phi$ such that $\vdash \phi$ is provable.  Some examples are shown in Figure~\ref{fig:theorems}.

\begin{figure}
  \centering
\begin{tabular}{lll} 
  $\gamma$ & & \text{Theorem over $\gamma$} \\ \hline
  $X$ & & If $X$ is a set, then $X$ is well-orderable. \\
        $x,y,z$ & & If  $x,y,z$ are ordinals, then $(x+y)+z = x+(y+z)$. \\
               $x$ & & $x=x$. 
\end{tabular}
\caption{Some theorems of Open \tops{}}
\label{fig:theorems}
\end{figure}
To make sense of \tops{} over $\gamma$ from a particularist viewpoint, let $\rho$ be a \emph{$\gamma$-valuation}, which associates to each $x \in \gamma$ a thing $\rho(x) \in \totall$.  Note that $\rho(x)$ may be an urelement, a Quine atom, an inaccessible cardinal, a measurable cardinal, etc.   With respect to $\rho$, every formula or sequent over $\gamma$ is bivalent, every instance of an axiom scheme is true, and every instance of a logical rule preserves truth; so every provable sequent and every theorem over $\gamma$ is true.

Accordingly, although $\totall$ is not fixed, and we are free to admit to it anything we find credible and desirable, we would not admit things that violate a theorem of Open \tops{}.  For example, we would not admit a set that is not well-orderable, ordinals on which addition is not associative, or a thing that is not equal to itself.  We deem such properties impossible, as they contradict the beliefs we have accepted.


\section{Related work} \label{sect:related}
\subsection{General background}

The following is just a selection, as the relevant literature is large.
\begin{itemize}
\item
For a wide-ranging account of set-theoretic belief, see Maddy~\cite{Maddy:believeI}.

\item For recent discussion of unrestricted quantification and the totality of sets, see e.g.\ the anthology \cite{RayoUzquiano:absgen} and the overviews~\cite{Florio:unrestricted,Studd:genextparadox}.  These issues have been considered in many contexts, such as reflection principles~\cite{HorstenWelch:reflectabsinf}, categoricity theorems~\cite{ButtonWalsh:structurecategoricity,Isaacson:realitymath,Martin:multiuniv}, modal logic~\cite{HamkinsLinnebo:modallogsetpot} and categorical semantics~\cite{Shulman:stacksemmatstruct,AwodeyButzSimpsonStreicher:relclasses}
  .  




\item   Following~\cite{Levy:hierarchy}, it is
   usual to classify sentences by their number of alternations of unrestricted quantifiers.   Mathias~\cite{Mathias:strengthmac} extensively analyzes subsystems of ZFC, such as the theories of Mac Lane~\cite{maclane:formandfunction} and Kripke-Platek, that use this classification to limit (in various ways) the permissible use of Separation and/or
   Replacement.  These are classical first-order theories, but other authors restrict also the use of excluded middle~\cite{Tharp:quasiint,Pozsgay:liberal,Pozsgay:semiint,Wolf:thesis,Friedman:someappkleene,Friedman:hanf,Feferman:strength,Rathjen:scopefeferman,Sanchis:settheoryoper}.  The legitimacy of classical vs intuitionstic reasoning for the totality of sets is further discussed in~\cite{Lear:setsem,Paseau:openend,LinneboShapiro:actupotinfty,Rumfitt:boundarythought}.  
  
\item  The formation of sequents in \tops{} is adapted from dependent type theory.  Indeed, extensional dependent type theory~\cite{MartinLoef:inttypetheory}, without universes, can be seen as a subsystem of \tops{}.

\item Set theories in which terms and propositional formulas are simultaneously defined can be seen e.g.\ in~\cite{Sazonov:boundedset,Moczydlowski:depset,Pozsgay:semiint}.

\item  The important role of well-founded recursion in set theory has been argued e.g.\ in~\cite{Wolf:efficienttransfrec,Mathias:strengthmac,Rin:transfiniterec}.

\item The convention that ``nonsense denotes the empty set'' is followed by the proof assistants Isabelle/ZF~\cite[Section 7.2]{Paulson:setverifI} and  Metamath~\cite{Metamath:iotanul}.
  
\end{itemize}

\subsection{Comparison with Mayberry's Local ZF}

After writing the first version of this paper, the author learned of the system ``Local ZF'' studied by Mayberry~\cite{Mayberry:globalquant,Mayberry:consistsetII}, with variants appearing in~\cite{Mayberry:newbegriffII,Mayberry:book}.  It clearly proves the same sentences (modulo notational differences) as \tops{}.  Therefore the following valuable results about Local ZF given in~\cite{Mayberry:globalquant} apply also to \tops{}.
\begin{itemize}
\item It has a conservative extension that provides  global choice function.
\item It has a further conservative extension that allows formulas with unrestricted quantification, governed by intuitionistic logic \emph{without induction} and various set-theoretic axioms. 
    \item It can be embedded in a weak fragment of \zf{}, whose consistency is provable in \zf{}.
   \end{itemize}
Despite the close similarity, Local ZF differs from \tops{} in several ways.  It adopts the purity assumptions, unlike \tops{}, and lacks the modular arrangement mentioned in Section~\ref{sect:subsyst}.  More importantly, it is not quite a theory of particular sets, as it uses \emph{open formulas} rather than sequents. 

To see the problem this causes, recall first the conventional understanding of open formula systems.  A formula's denotation is given with respect to a \emph{valuation}, i.e.\ map $\vars \to \totall$.  A \emph{universally true} formula is one that is true with respect to every valuation.  A theory is deemed acceptable when each axiom is universally true and each inference rule preserves universal truth.

For example, consider the following Hilbert-style counterpart of our $\forall_{\in}\Intro$ rule:
\begin{displaymath}
  \begin{prooftree}
    \psi \Rightarrow (x \smin A \Rightarrow \phi(x))
    \using x \text{ fresh for } \psi,A,\phi
    \justifies
    \psi \Rightarrow \forall x \smin A.\,\phi(x)
  \end{prooftree}
\end{displaymath}
(Rule III(3) in~\cite{Mayberry:globalquant} is similar.)   For a totalist, it clearly preserves universal truth.  But for a particularist, who doubts the notion of universal truth, it is hard to see what kind of truth this rule can be said to preserve.

  By contrast, the logical rules for \tops{} preserve not \emph{universal truth} but simply \emph{truth}---that is what makes them immediately convincing.  Sequents are bivalent just as sentences are, because every variable is declared.  Specifically, the premise of the $\forall_{\in}\Intro$ rule in Figure~\ref{fig:logrules} would not be a sequent if $x \colin A$ were replaced by $x \smin A$.  

   What about Open \tops{}? The reason \tops{} over $\gamma$ provides an acceptable way of reasoning about a particular $\gamma$-valuation $\rho$ is that each axiom is true with respect to $\rho$, and each rule preserves truth with respect to $\rho$.  No other valuation need be considered.

   Thus, whereas one can present first-order logic either in Hilbert style or via Natural Deduction, this is not the case for \tops{}.  The latter must be presented via Natural Deduction, distinguishing between a declaration $x\colin A$ and a mere hypothesis of the form $x \smin A$, in order to have a notion of truth that every rule preserves.

\section{Conclusion} \label{sect:conclude}

Despite the popularity of \zf{} as a set-theoretic foundation, the Burali-Forti paradox (for example) raises serious concerns about its meaningfulness.  As an alternative, we propose \tops{}, a Natural Deduction system that avoids the problem by speaking only about particular sets.  This allows us to maintain, as particularists, that its sentences are bivalent and its theorems objectively true.  An example theorem is the well-orderability of $V_{\omega_1}$.  Theorems of Open \tops{} have a different character; they lack objective meaning but are true with respect to any valuation.

\noindent \paragraph*{Acknowledgements} I thank Ehud Hrushovski, Adrian Mathias and Norman Megill for helpful comments.


\begin{thebibliography}{10}

\bibitem{AntosFriedmanHonzikTernullo:hyperun}
C.~Antos, S.-D. Friedman, R.~Honzik, and C.~Ternullo, editors.
\newblock {\em The Hyperuniverse Project and Maximality}.
\newblock Birkh\"{a}user, 2018.

\bibitem{AwodeyButzSimpsonStreicher:relclasses}
Steven Awodey, Carsten Butz, Alex Simpson, and Thomas Streicher.
\newblock Relating first-order set theories, toposes and categories of classes.
\newblock {\em Ann. Pure Appl. Logic}, 165(2):428--502, 2014.

\bibitem{Bernays:platonism}
P.~Bernays.
\newblock Platonism in mathematics, {English} translation of a 1934 lecture.
\newblock In P.~Benacerraf and H.~Putnam, editors, {\em Philosophy of
  mathematics: selected readings}. Prentice-Hall, 1964.

\bibitem{Boolos:whencecontra}
George Boolos.
\newblock Whence the contradiction?
\newblock {\em Aristotelian Society Supplementary Volume}, 67:211--233, 1993.

\bibitem{ButtonWalsh:structurecategoricity}
Tim Button and Sean Walsh.
\newblock Structure and categoricity: Determinacy of reference and truth value
  in the philosophy of mathematics †.
\newblock {\em Philosophia Mathematica}, 24(3):283--307, 2016.

\bibitem{Cartwright:speakingevery}
Richard~L. Cartwright.
\newblock Speaking of everything.
\newblock {\em No\^{u}s}, 28(1):1--20, 1994.

\bibitem{Dummett:fregephilmath}
M.~Dummett.
\newblock {\em Frege: Philosophy of Mathematics}.
\newblock Harvard University Press, 1991.

\bibitem{Feferman:strength}
S.~Feferman.
\newblock On the strength of some semi-constructive theories.
\newblock In U.~Berger, P.~Schuster, and M.~Seisenberger, editors, {\em Logic,
  Construction, Computation}. Ontos, 2012.

\bibitem{Feferman:predicativism}
Solomon Feferman.
\newblock Predicativity.
\newblock In Stewart Shapiro, editor, {\em Oxford Handbook of Philosophy of
  Mathematics and Logic}, pages 590--624. Oxford: Oxford University Press,
  2005.

\bibitem{Florio:unrestricted}
Salvatore Florio.
\newblock Unrestricted quantification.
\newblock {\em Philosophy Compass}, 9(7):441--454, 2014.

\bibitem{Friedman:someappkleene}
H.~Friedman.
\newblock Some applications of {Kleene's} methods for intuitionistic systems.
\newblock In A.R.D.Mathias and H.Rogers, editors, {\em Cambridge Summer School
  in Mathematical Logic 1971}, volume 337 of {\em Lecture Notes in
  Mathematics}. Springer, 1973.

\bibitem{Friedman:hanf}
Harvey Friedman.
\newblock On existence proofs of {Hanf} numbers.
\newblock {\em J. Symb. Log.}, 39(2):318--324, 1974.

\bibitem{Friedman:evidencehyperun}
Sy-David Friedman.
\newblock Evidence for set-theoretic truth and the hyperuniverse programme.
\newblock In Antos et~al. \cite{AntosFriedmanHonzikTernullo:hyperun}, pages
  75--107.

\bibitem{HamkinsLinnebo:modallogsetpot}
Joel~David Hamkins and \O{}ystein Linnebo.
\newblock The modal logic of set-theoretic potentialism and the potentialist
  maximality principles.
\newblock {\em to appear in Review of Symbolic Logic}, 2019.

\bibitem{Isaacson:realitymath}
Daniel Isaacson.
\newblock The reality of mathematics and the case of set theory.
\newblock In Z.~Novak and A.~Simonyi, editors, {\em Truth, Reference and
  Realism}, pages 1--76. Central European Press, Budapest, 2011.

\bibitem{Lear:setsem}
Jonathan Lear.
\newblock Sets and semantics.
\newblock {\em The Journal of Philosophy}, 74(2):86--102, 1977.

\bibitem{Levy:hierarchy}
A.~{Levy}.
\newblock {A hierarchy of formulas in set theory.}
\newblock {\em {Mem. Am. Math. Soc.}}, 57:76, 1965.

\bibitem{LinneboShapiro:actupotinfty}
{\O}ystein Linnebo and Stewart Shapiro.
\newblock Actual and potential infinity.
\newblock {\em No\^{u}s}, 53(1):160--191, 2019.

\bibitem{maclane:formandfunction}
Saunders Mac~Lane.
\newblock {\em Mathematics, form and function}.
\newblock Springer-Verlag, New York, 1986.

\bibitem{Maddy:believeI}
P.~Maddy.
\newblock Believing the axioms {I}.
\newblock {\em Journal of Symbolic Logic}, 53(2):481--511, 1988.

\bibitem{Martin:multiuniv}
Donald~A. Martin.
\newblock Multiple universes of sets and indeterminate truth values.
\newblock {\em Topoi}, 20(1):5--16, 2001.

\bibitem{MartinLoef:inttypetheory}
P.~Martin-{L}{\"o}f.
\newblock {\em Intuitionistic type theory}.
\newblock Bibliopolis, Napoli, 1980.

\bibitem{Mathias:strengthmac}
A.~R.~D. Mathias.
\newblock The strength of {Mac Lane} set theory.
\newblock {\em Ann. Pure Appl. Logic}, 110(1-3):107--234, 2001.

\bibitem{Mayberry:book}
J.~P. Mayberry.
\newblock {\em The Foundations of Mathematics in the Theory of Sets}.
\newblock Cambridge University Press, 2000.

\bibitem{Mayberry:consistsetII}
John Mayberry.
\newblock The consistency problem for set theory: An essay on the {Cantorian}
  foundations of mathematics ({II}).
\newblock {\em British Journal for the Philosophy of Science}, 28(2):137--170,
  1977.

\bibitem{Mayberry:globalquant}
John Mayberry.
\newblock Global quantification in {Zermelo-Fraenkel} set theory.
\newblock {\em Journal of Symbolic Logic}, 50(2):289--301, 1985.

\bibitem{Mayberry:newbegriffII}
John~P. Mayberry.
\newblock A new {Begriffsschrift (II)}.
\newblock {\em British Journal for the Philosophy of Science}, 31(4):329--358,
  1980.

\bibitem{McNaughton:conceptualschemes}
Robert McNaughton.
\newblock Conceptual schemes in set theory.
\newblock {\em The Philosophical Review}, 66(1):66--80, 1957.

\bibitem{Moczydlowski:depset}
Wojciech Moczydlowski.
\newblock A dependent set theory.
\newblock In {\em 22nd Annual {IEEE} Symposium on Logic in Computer Science},
  pages 23--34. IEEE Computer Society, 2007.

\bibitem{Parsons:setsandclasses}
Charles Parsons.
\newblock Sets and classes.
\newblock {\em No\^{u}s}, 8(1):1--12, 1974.

\bibitem{Paseau:openend}
A.~Paseau.
\newblock The open-endedness of the set concept and the semantics of set
  theory.
\newblock {\em Synthese}, 135(3):379--399, 2003.

\bibitem{Paulson:setverifI}
Lawrence~C. Paulson.
\newblock Set theory for verification: {I.} {F}rom foundations to functions.
\newblock {\em Journal of Automated Reasoning}, 11(3):353--389, December 1993.

\bibitem{PelletierHazen:histnatded}
Francis~Jeffry Pelletier and Allen~P. Hazen.
\newblock A history of natural deduction.
\newblock In Dov~M. Gabbay, Francis~Jeffry Pelletier, and John Woods, editors,
  {\em Logic: {A} History of its Central Concepts}, volume~11 of {\em Handbook
  of the History of Logic}, pages 341--414. Elsevier, 2012.

\bibitem{Pozsgay:semiint}
Lawrence~J. Pozsgay.
\newblock Semi-intuitionistic set theory.
\newblock {\em Notre Dame Journal of Formal Logic}, 13(4):546--550, 1972.

\bibitem{Pozsgay:liberal}
L.J. Pozsgay.
\newblock Liberal intuitionism as a basis for set theory.
\newblock In D.S.Scott, editor, {\em Proceedings, Axiomatic Set Theory I (Part
  1)}, volume XIII of {\em Proceedings of Symposia in Pure Mathematics}, pages
  321--330. American Mathematical Society, 1971.

\bibitem{Rathjen:scopefeferman}
Michael Rathjen.
\newblock The scope of {Feferman's} semi-intuitionistic set theories and his
  second conjecture.
\newblock {\em Indagationes Mathematicae}, 30(3):500 -- 525, 2019.

\bibitem{RayoUzquiano:absgen}
A.~Rayo and G.~Uzquiano, editors.
\newblock {\em Absolute generality}.
\newblock Oxford University Press, 2006.

\bibitem{Rin:transfiniterec}
Benjamin~G. Rin.
\newblock Transfinite recursion and computation in the iterative conception of
  set.
\newblock {\em Synthese}, 192(8):2437--2462, 2015.

\bibitem{Rumfitt:boundarythought}
Ian Rumfitt.
\newblock {\em The Boundary Stones of Thought: An Essay in the Philosophy of
  Logic}.
\newblock Oxford University Press, United Kingdom, 2015.

\bibitem{Metamath:iotanul}
A.~Salmon, 2011.
\newblock \texttt{http://us.metamath.org/mpeuni/iotanul.html}.

\bibitem{Sanchis:settheoryoper}
Luis~E. Sanchis.
\newblock {\em Set Theory - An Operational Approach}.
\newblock Gordon and Breach Science Publishers, 1996.

\bibitem{Sazonov:boundedset}
V.Yu. Sazonov.
\newblock On bounded set theory, invited talk.
\newblock In {\em 10th International Congress on Logic, Methodology and
  Philosophy of Sciences in Volume I: Logic and Scientific Method}, pages
  85--103. Kluwer, 1997.

\bibitem{Shulman:stacksemmatstruct}
M.~Shulman.
\newblock Stack semantics and the comparison of material and structural set
  theories.
\newblock Available at {\tt arxiv.org/abs/1004.3802}, 2010.

\bibitem{Studd:genextparadox}
J.P. Studd.
\newblock {V}---generality, extensibility, and paradox.
\newblock {\em Proceedings of the Aristotelian Society}, 117(1):81--101, 2017.

\bibitem{Tharp:quasiint}
Leslie~H. Tharp.
\newblock A quasi-intuitionistic set theory.
\newblock {\em J. Symb. Log.}, 36(3):456--460, 1971.

\bibitem{Velleman:constructlib}
Daniel~J. Velleman.
\newblock {Constructivism Liberalized}.
\newblock {\em Philosophical Review}, 102(1), 1993.

\bibitem{HorstenWelch:reflectabsinf}
Philip Welch and Leon Horsten.
\newblock Reflecting on absolute infinity.
\newblock {\em Journal of Philosophy}, 113(2):89--111, 2016.

\bibitem{Wolf:thesis}
R.~S. Wolf.
\newblock {\em Formally intuitionistic set theories with bounded predicates
  decidable}.
\newblock PhD thesis, Stanford University, 1974.

\bibitem{Wolf:efficienttransfrec}
Robert~S. Wolf.
\newblock A highly efficient ``transfinite recursive definitions'' axiom for
  set theory.
\newblock {\em Notre Dame Journal of Formal Logic}, 22(1):63--75, 1981.

\end{thebibliography}

\end{document}